\documentclass{elsart}

\usepackage{amssymb}

\newcommand{\dT}{d_{\mathrm{T}}}
\newcommand{\id}{\mathrm{id}}
\renewcommand{\Im}{\mathop{\mathrm{Im}}}
\newcommand{\pd}[2]{\frac{\partial #1}{\partial #2}}
\newcommand{\pdb}[3]{\frac{\partial #1}{\partial #2 \, \partial #3}}

\newcommand{\quarter}{{\textstyle\frac{1}{4}}}
\newcommand{\R}{\mathbb{R}}
\newcommand{\To}{T_\circ}
\newcommand{\vf}[1]{\frac{\partial}{\partial #1}}

\begin{document}

\begin{frontmatter}

\title{How to recover a Lagrangian using the homogeneous variational bicomplex}

\author{D.J.\ Saunders\corauthref{cor1}\thanksref{label1}}
\corauth[cor1]{Address for correspondence: 30 Little Horwood Road, Great Horwood, Milton Keynes, MK17 0QE, UK}
\thanks[label1]{This work was supported by grant no. 201/06/0922 for Global Analysis and its Applications from the Czech Science Foundation.}

\address{Department of Algebra and Geometry,\\ Palack\'{y} University, Olomouc, Czech Republic}

\ead{david@symplectic.demon.co.uk}

\begin{abstract}
We show how the homogeneous variational bicomplex provides a useful formalism for describing a number of properties of single-integral variational problems, and we introduce a subsequence of one of the rows of the bicomplex which is locally exact with respect to the variational derivative. We are therefore able to recover a Lagrangian from a set of equations given as a variationally-closed differential form. As an example, we show how to recover a first-order Lagrangian from a suitable set of second-order equations.
\end{abstract}

\begin{keyword}
variational bicomplex

\MSC 58E99 \sep 49F99
\end{keyword}
\end{frontmatter}

\section{Introduction}

I met Willy~Sarlet in 1986 in Gent, at the first of the annual Workshops on Differential Geometric Methods in Mechanics; since then we have had many discussions on topics of mutual interest, and worked together (along with Frans~Cantrijn) on a series of papers on non-holonomic mechanics. The present paper is on a rather different subject, one which I have been studying on and off for some time. In fact I did not give a talk at the meeting to honour Willy's sixtieth birthday, as he had already heard parts of this paper twice before and I did not feel that it would be appropriate to trouble him with it again. But at the meeting itself there was some discussion of the equations to be satisfied when `Euler-Lagrange forms' arise from a Lagrangian in the fibred case, and so I was prompted to see how these equations might work in the homogeneous context. These new ideas are based on the constructions I had described previously, and can be expressed in a manner which is perhaps more elegant than in the fibred case. I shall restrict attention here to the case of a single-integral variational problem, although much of the discussion may be extended, albeit in a significantly more complicated manner, to multiple-integral problems.

My interest in homogeneous variational problems was encouraged by Mike Crampin, and some aspects of this work, particularly that described in Section~\ref{sec:hom}, may be found in three joint papers published a few years ago~\cite{CS1,CS2,CS3}; these papers also cover multiple-integral problems.

\section{Homogeneous variational problems}
\label{sec:hom}

There is an argument for saying that variational problems fall into two distinct types, depending upon the nature of the solutions. In one type, the solutions will be parametrised submanifolds of the ambient manifold; in the other, they will be submanifolds without any preferred parametrisation.

Different geometrical structures are used to describe these two types of problem. In the first case, we would typically have a fibration $\pi : E \to M$ and consider a Lagrangian on a jet manifold $J^k\pi$; the solutions are local sections $\phi$ of $E \to M$, and the solution submanifolds $\Im(\phi)$ are parametrised by the sections themselves. Problems of this type arise in mathematical physics:\ for instance, in classical mechanics, where typically $M = \R$, and in field theories, where $\dim M > 1$.

The other type of problem is that of minimal submanifolds, and in the first-order single-integral case this is just Finsler geometry. Here the ambient manifold $E$ has no fibration, and instead of the affine jet manifolds $J^k\pi$ we could use instead the manifolds of $m$-dimensional contact elements $J^k(E, m)$ (otherwise known as jets of immersions, or as higher-order Grassmannians); the particular case $J^1(E, 1)$ is just the projective tangent bundle of $E$. But another approach, described in our papers~\cite{CS1,CS2,CS3} and elsewhere, is to use instead the manifolds of non-degenerate $m$-velocities $T^k_{\circ m}E$. These are the total spaces of principal bundles $T^k_{\circ m}E \to J^k(E, m)$ with structure groups $L^k_m$, the $k$-th order $m$-dimensional jet groups; using these manifolds means that we should consider homogeneous objects, with a specific behaviour under the group action. There is an argument that working directly on $J^k(E, m)$ is more straightforward, as there is no need to check invariance. But actually that is a disadvantage, as starting with homogeneous objects one often finds that the resulting constructions are \emph{not} invariant under the group action, and indeed this fact explains some of the difficulties found in the case of higher-order field theories.

There is, of course, a relationship between the two types of problem. The `homogenisation trick' has been known for many years in classical mechanics:\ here, one takes a Lagrangian $L(t, q^i, \dot{q}^i)$ given as a function with respect to the standard volume form $dt$, and replaces it by a homogeneous Lagrangian
\[
\widetilde{L}(t, \dot{t}, q^i, \dot{q}^i) = \dot{t} L(t, q^i, \dot{q}^i \dot{t}^{-1})  
\]
where $t$ is now regarded as a variable with the same status as the $q^i$. Geometrically, one has moved from the jet manifold $J^1\pi$ (where $\pi : E \to \R$ is some fibred manifold) to a suitable open subset of the slit tangent manifold $\To E$, essentially by taking the contraction of $L \, dt$ with the total time derivative. Much the same sort of trick works for higher orders, and also for multiple integral problems.

In this paper we shall concentrate on the case $m=1$, and so consider the higher-order non-degenerate tangent bundles $\To^k E$. We take as known (see, for instance,~\cite{CSC}) the existence of two canonical constructions on these manifolds. The first of these is the \emph{total time derivative $\dT$}, a vector field along the projection $\To^{k+1} E \to \To^k E$ expressed in coordinates as
\[
\dT = \sum_{p=0}^k q^i_{(p+1)} \vf{q^i_{(p)}}
\]
where $q^i$ are local coordinates on $E$ and $q^i_{(p)}$ are the corresponding derivative coordinates with `$p$ dots' on a tangent manifold. The second is the \emph{vertical endomorphism $S$}, a type $(1,1)$ tensor field on $\To^k E$ expressed in coordinates as
\[
S = \sum_{p=0}^{k-1} (p+1) dq^i_{(p)} \otimes \vf{q^i_{(p+1)}} \, .
\]

In order to express conditions of homogeneity, we need to make use of some more vector fields $\Delta^p$ on $\To^k E$; these are defined as the contractions $\Delta^p = S^p (\dT)$ for $1 \leq p \leq k$. Note that, although $\dT$ is a vector field along a projection, the properties of $S$ ensure that $\Delta^p$ are genuine vector fields on a manifold:\ they are, indeed, the fundamental vector fields of the principal bundle $\To^k E \to J^k(E,1)$.

A \emph{Lagrangian} in this context is a function on $\To^k E$. It is not a differential form:\ rather, it defines the variational problem
\[
\delta \int (j^k \gamma)^*(L) dt = 0
\]
where $\gamma$ is a curve in $E$, $j^k \gamma$ is its prolongation to a curve in $\To^k E$, and the integral is taken with respect to the standard volume form $dt$ on $\R$. The Lagrangian is \emph{homogeneous} if $\Delta^1(L) = L$, $\Delta^p(L) = 0$ for $2 \leq p \leq k$.

Each Lagrangian gives rise to its \emph{Hilbert form}
\[
\vartheta_L = \sum_{p=0}^{k-1} \frac{(-1)^p}{(p+1)!} \dT^p S^{p+1} dL
\]
defined on $\To^{2k-1} E$; this 1-form has the same extremals as $L$~\cite{CS2} and is a generalisation of the Hilbert form used in Finsler geometry, where $k=1$ and the formula simplifies to $\vartheta_L = SdL$. If $L$ is a homogeneous Lagrangian then $\vartheta_L$ is projectable to $J^{2k-1}(E,1)$. If $L$ (perhaps defined on an open submanifold of $\To^k E$) has been obtained from a Lagrangian form on some $J^k\pi$ by means of the homogenisation trick then the projected image of $\vartheta_L$ is just the Poincar\'{e}-Cartan form of the original Lagrangian.

Finally, the Lagrangian and its Hilbert form may be used together to construct the Euler-Lagrange form $\varepsilon_L$ using the formula
\[
\varepsilon_L = dL - \dT \vartheta_L \, ;
\]
we therefore see, explicitly, that
\[
\varepsilon_L = \sum_{p = 0}^k \frac{(-1)^p}{p!} \dT^p S^p dL \, .
\]

\section{The homogeneous variational bicomplex}

The cohomological approach to the study of variational problems involves constructing suitable sequences of differential forms, or alternatively a bicomplex of differential forms. Historically the bicomplex appeared first, and involved forms defined on an infinite jet manifold~\cite{AndDuch,Tak,Tul,Vin1,Vin2}. The use of a single sequence appeared a little later, and used forms defined on a finite-order jet manifold~\cite{KrupSeq}. Subsequently, versions of the bicomplex on finite-order jet manifolds were also obtained~\cite{Vit}. All these constructs were defined initially for jets of local sections of a fibration; versions involving contact elements have also appeared, and in some cases the elements of the spaces have been equivalence classes of forms rather than individual forms.

The version of the bicomplex suitable for homogeneous problems was described in~\cite{Sau}, and in its full generality uses certain spaces of vector-valued forms rather than scalar forms. But for single-integral variational problems the target vector spaces are all one-dimensional, and so we may consider the vector-valued forms as identified with scalar forms; that is the approach adopted here. The manifolds on which the forms are defined are, as before, the manifolds of non-degenerate higher-order tangent vectors $\To^k E$ on some given manifold $E$ with $\dim E = n+1$; we write $\Omega^r_k$ for the space of $r$-forms on $\To^k E$. (Later on we shall talk about certain operators being locally exact, and then it might be more appropriate to regard $\Omega^r_k$ as the sheaf of germs of $r$-forms; we leave the reader to make this interpretation as required.)

The basic idea is to combine the usual de Rham sequence of forms with the total time derivative operator $\dT$, as in Figure 1. (It is convenient here to modify the start of the sequence, replacing the usual $0 \to \R \to \Omega^0$ by $0 \to \overline{\Omega}^0$, where $\overline{\Omega}^0 = \Omega^0 / \R$.) In this diagram we have been explicit about the order of the forms; but in fact the operator used to demonstrate exactness of the columns will, in general, increase the order:\ indeed, `exactness' there must be understood as being given modulo a pull-back to a higher-order manifold, and so we should imagine not a single bicomplex but a whole family, indexed by order. We shall therefore tend to omit any reference to order in the main body of the paper, and so improve the clarity of the formul\ae. In any particular case, one may readily calculate the order (or, more accurately, an upper bound on the order) of the forms concerned; the general principle is to `keep going until you run out of dots'.

\begin{figure}
\begin{center}
\begin{picture}(300,350)(50,-20)
\multiput(100,300)(100,0){3}{\vector(0,-1){40}}
\multiput(100,220)(100,0){3}{\vector(0,-1){40}}
\multiput(100,140)(100,0){3}{\vector(0,-1){40}}
\multiput(100,60)(100,0){3}{\vector(0,-1){40}}
\multiput(40,240)(0,-80){2}{\vector(1,0){35}}
\put(40,80){\vector(1,0){20}}
\multiput(130,240)(0,-80){3}{\vector(1,0){35}}
\multiput(230,240)(0,-80){3}{\vector(1,0){35}}
\multiput(330,240)(0,-80){3}{\vector(1,0){35} \quad $\ldots$ }
\multiput(20,240)(0,-80){3}{\makebox(0,0){$0$}}
\multiput(100,320)(100,0){3}{\makebox(0,0){$0$}}
\put(100,240){\makebox(0,0){$\overline{\Omega}^{\;0}_k$}}
\put(200,240){\makebox(0,0){$\Omega^1_k$}}
\put(300,240){\makebox(0,0){$\Omega^2_k$}}
\put(100,160){\makebox(0,0){$\overline{\Omega}^{\;0}_{k+1}$}}
\put(200,160){\makebox(0,0){$\Omega^1_{k+1}$}}
\put(300,160){\makebox(0,0){$\Omega^2_{k+1}$}}
\put(100,80){\makebox(0,0){$\overline{\Omega}^{\;0}_{k+1} / \dT \overline{\Omega}^{\;0}_k$}}
\put(200,80){\makebox(0,0){$\Omega^1_{k+1} / \dT \Omega^1_k$}}
\put(300,80){\makebox(0,0){$\Omega^2_{k+1} / \dT \Omega^2_k$}}
\multiput(100,0)(100,0){3}{\makebox(0,0){$0$}}
\multiput(150,170)(0,80){2}{\makebox(0,0)[b]{$\scriptstyle d$}}
\multiput(250,170)(0,80){2}{\makebox(0,0)[b]{$\scriptstyle d$}}
\multiput(105,200)(100,0){3}{\makebox(0,0)[l]{$\scriptstyle\dT$}}
\end{picture}
\end{center}
\caption{The homogeneous variational bicomplex}
\end{figure}

Our first task is to investigate the exactness of this bicomplex. Define the differential operator $P : \Omega^r \to \Omega^r$ (with $r \geq 1$) by
\[
P = \sum_p \frac{(-1)^p}{r^{p+1} (p+1)!} \dT^p S^{p+1} \, ,
\]
where the sum ranges from zero to a suitably large number so that $S^{p+1}\theta$ vanishes:\ if one wished to be explicit one would take a sum up to $rk-1$.
\begin{thm}
For every $r$-form $\theta \in \Omega^r$
\[
P \dT \theta = \theta \quad \mbox{modulo a pull-back.}
\]
\end{thm}
\begin{pf}
We use the commutation relation
\[
S\dT\theta = \dT S\theta + r\theta \qquad \mbox{modulo a pull-back}
\]
which may be established easily in coordinates for a 1-form and then extended to an $r$-form by induction. Thus
\[
S^{p+1}\dT\theta = \dT S^{p+1}\theta + r(p+1)S^p\theta \, ,
\]
and so
\[
\dT^p S^{p+1}\dT\theta = \dT^{p+1}S^{p+1}\theta + r(p+1)\dT^p S^p\theta \, .
\]
The coefficients in the alternating sum are chosen so that successive terms will cancel out, leaving just the second part of the term for $p = 0$ which gives $\theta$ as required. \qed
\end{pf}

It follows from this that the columns of the bicomplex, apart perhaps from the first, are globally exact (modulo pull-backs, of course). The first two rows are locally exact as they are the usual de Rham sequences; from this we may show by chasing round diagrams that the first column and the bottom row are also locally exact.  

Now the bottom row of the bicomplex contains equivalence classes of forms, differing by total time derivatives. It is convenient to be able to identify a representative form in each class, and the differential operator $P$ allows us to do this, at the expense of a possible increase in the order of the forms.
\begin{lem}
\label{Lem:canrep}
For any $r$-form $\theta \in \Omega^r$ the $r$-form
\[
\theta - \dT P\theta
\]
is a canonical representative of the class $[\theta] \in \Omega^r / \dT\Omega^r$.
\end{lem}
\begin{pf}
It is clear that $\theta - \dT P\theta \in [\theta]$, so we just have to show that starting with a different form $\hat{\theta} \in [\theta]$ we would achieve the same result. But by definition
\[
\hat{\theta} = \theta + \dT\phi
\]
for some $\phi \in \Omega^r$, and 
\[
(\dT \phi) - \dT P (\dT \phi) = 0
\]
as $P\dT\phi = \phi$. \qed
\end{pf}

It is now possible to give a very simple description of the canonical representative.
\begin{thm}
\label{thm:Psi-space}
Let $\Psi^r \subset \Omega^r$ (with $r \geq 1$) be the subspace of $r$-forms defined by
\[
\Psi^r = \{ \theta \in \Omega^r : S\theta = 0 \} \, .
\]
Then $\id - \dT P$ is a projection operator $\Omega^r \to \Psi^r$, so that $\theta \in \Omega^r$ is the canonical representative of its class $[\theta]$ if, and only if, $\theta \in \Psi^r$.
\end{thm}
\begin{pf}
It is an immediate consequence of Lemma~\ref{Lem:canrep} that $\id - \dT P$ is a projection operator on $\Omega^r$, so we need to confirm that its image is $\Psi^r$. From the formula for $P$ we see that if $\theta \in \Psi^r$ then $P\theta = 0$, so that $\theta = (\id - \dT P)\theta$ and therefore that $\Psi^r \subset \Im(\id - \dT P)$.

For the converse, suppose that $\theta = (\id - \dT P) \phi$ for some $\phi \in \Omega^{r-1}$. We note first that
\[
S \dT^p = \dT^{p+1}S + r(p + 1)\dT^p
\]
from the commutation relation mentioned earlier, so that
\begin{eqnarray*}
S\dT P\phi & = & S\dT \sum_p \frac{(-1)^p}{r^{p+1}(p+1)!} \dT^p S^{p+1} \phi \\
& = & \sum_p \frac{(-1)^p}{r^{p+1}(p+1)!} (\dT^{p+1} S^{p+2}+ r(p + 1)\dT^p S^{p+1})\phi \\
& = & S\phi
\end{eqnarray*}
as the alternating sum collapses once again. We conclude that
\[
S\theta = S(\phi - \dT P\phi) = S\phi - S\phi = 0
\]
so that $\theta \in \Psi^r$ as required.
\qed
\end{pf}
We note that, in the case $r=1$, $S\theta = 0$ exactly when $\theta$ is horizontal over the manifold $E$. This is not, of course, true when $r \geq 2$.

\section{The variational derivative}

The properties of the operator $\id - \dT P$ suggest that we should consider replacing the de Rham differential $d : \Omega^r \to \Omega^{r+1}$ by a new operator
\[
\delta : \Omega^r \to \Omega^{r+1} \, , \qquad \delta\theta = d\theta - \dT Pd\theta 
\]
which we shall call the \emph{variational derivative}. The new operator $\delta$ is clearly $\R$-linear, although --- despite its name --- it is not in fact a derivation. It is, however, a coboundary operator.
\begin{lem}
The operator $\delta$ satisfies $\delta^2 = 0$.
\end{lem}
\begin{pf}
From $d^2 = 0$ and $P\dT = \id$ we obtain
\begin{eqnarray*}
(d - \dT P d)(d - \dT P d) & = & - d\dT P d + \dT P d \dT P d \\
& = &  - \dT d P d + \dT P \dT d P d \\
& = & 0 \, .
\end{eqnarray*}
\end{pf}
The operator $\delta$ mapping $\Omega^r$ to $\Omega^{r+1}$ is not locally exact for any $r$, even modulo pullbacks. But by construction $\delta$ takes its values in $\Psi^{r+1} \subset \Omega^{r+1}$, so we may consider the restriction of $\delta$ to $\Psi^r$:
\[
\overline{\Omega}^0 \stackrel{\delta}{\to} \Psi^1 \stackrel{\delta}{\to} \Psi^2 
\stackrel{\delta}{\to} \cdots \stackrel{\delta}{\to} \Psi^r \stackrel{\delta}{\to}
\Psi^{r+1} \stackrel{\delta}{\to} \cdots \stackrel{\delta}{\to} 0 \; .
\] 
This new sequence is indeed locally exact. To prove this, we use the homotopy formula for the de Rham differential
\[
h d\theta + dh\theta = \theta
\]
where $h$ is the usual Poincar\'{e} operator.
\begin{thm}
\label{thm:localexact}
If the $r$-form $\theta \in \Psi^r$ (with $r \geq 1$) satisfies $\delta\theta = 0$ then there is an $(r-1)$-form
\[
\psi \in \left\{
\begin{array}{ll}
\Psi^{r-1} & (r \geq 2) \\
\overline{\Omega}^0 & (r = 1)
\end{array}
\right.
\]
satisfying $\delta\psi = \theta$.
\end{thm}
\begin{pf}
From $\delta\theta = 0$ we have
\[
d\theta = \dT P d\theta \, ,
\]
so that
\[
\dT dPd\theta = d\dT Pd\theta = d^2\theta = 0 
\]
and therefore
\[
dPd\theta = P\dT dPd\theta = 0 \, .
\]
So $Pd\theta$ is $d$-closed, and locally we may put $\kappa = hPd\theta$ so that
\[
d\kappa = Pd\theta \, .
\]
We now have
\[
d\theta = \dT d\kappa = d\dT\kappa \, ,
\]
so that $\theta - \dT\kappa$ is $d$-closed and locally we may put $\phi = h(\theta - \dT\kappa) \in \Omega^{r-1}$ so that
\[
d\phi = \theta - \dT\kappa \, .
\]
Then
\begin{eqnarray*}
\delta\phi & = & d\phi - \dT P d\phi \\
& = & (\theta - \dT\kappa) - \dT P (\theta - \dT\kappa) \\
& = & \theta - \dT\kappa - \dT P \theta + \dT P\dT\kappa \\
& = & \theta - \dT P \theta \\
& = & \theta
\end{eqnarray*}
where the final equality arises because $\theta \in \Psi^r$ is its own canonical representative.

We now consider separately the cases $r = 1$ and $r \geq 2$. If $\theta \in \Omega^1$, so that $\phi \in \Omega^0$, we simply take $\psi = [\phi] \in \overline{\Omega}^0$, and then immediately $\delta\psi = \theta$. If instead $r \geq 2$ then $\phi \in \Omega^{r-1}$, so we put $\psi = \phi - \dT P \phi \in \Psi^{r-1}$, and we see that
\begin{eqnarray*}
\delta\psi & = & \delta (\phi - \dT P \phi) \\
& = & (\id - \dT P) d (\phi - \dT P \phi) \\
& = & d\phi - \dT Pd\phi - d\dT P\phi + \dT Pd\dT P\phi \\
& = & d\phi - \dT Pd\phi - \dT dP\phi + \dT P\dT dP\phi \\
& = & d\phi - \dT Pd\phi \\
& = & \theta
\end{eqnarray*}
as required. \qed
\end{pf}

\section{What has all this to do with the calculus of variations?}

The last two sections have been quite abstract, and we need to return to more concrete questions. We do this by relating the first three columns of the homogeneous variational bicomplex to specific aspects of the calculus of variations.

So consider a variational problem $\delta (\int L \, dt) = 0$, for some Lagrangian function $L \in \To^k$. (The Lagrangian might indeed be homogeneous, but for the purposes of the present discussion we do not require this.) We then have $[L] \in \overline{\Omega}^0$, and we may consider the variational derivative operator defined above. We have
\begin{eqnarray*}
\delta[L] & = & dL - \dT P dL \\
& = & dL - \dT \sum_p \frac{(-1)^p}{(p+1)!} \dT^p S^{p+1} dL \\
& = & \sum_{p^\prime} \frac{(-1)^{p^\prime}}{p^\prime!} \dT^{p^\prime} S^{p^\prime} dL
\end{eqnarray*}
where we have relabelled the sum with $p^\prime = p+1$; thus $\delta[L] = \varepsilon_L$, the Euler-Lagrange form for $L$ described in Section~\ref{sec:hom}. It immediate from Theorem~\ref{thm:Psi-space} that $S\varepsilon_L = 0$, so that $\varepsilon_L$ is horizontal over $E$.

We can also use the local exactness of $\delta$ to go backwards:\ given a 1-form $\varepsilon \in \Psi^1$, when does this represent the Euler-Lagrange form of some Lagrangian? This is the `simple' version of the Inverse Problem of the Calculus of Variations, where the multiplier matrix is fixed, and (as in the fibred case) there is a definite answer. We consider the \emph{Helmholtz-Sonin form} $\delta\varepsilon \in \Psi^2$, and if this vanishes then we will be able to construct a (local) Lagrangian for $\varepsilon$ using the formula from Theorem~\ref{thm:localexact}, which we may write explicitly as
\[
L = h(\varepsilon - \dT hPd\varepsilon) \, .
\]
The question of whether this Lagrangian has minimal order will depend on the choice of form $\kappa$ satisfying $d\kappa = Pd\varepsilon$.

As an example, we consider the particular case of a second-order form $\varepsilon \in \Psi^1_2$, where in coordinates we have $\varepsilon = \varepsilon_i dq^i$; we may think of the form $\varepsilon$ as a second-order differential operator taking its values in the vector bundle $T^* E$, and the differential equation itself as the kernel of the operator, a submanifold of $\To^2 E$. The coordinate formula for the Helmholtz-Sonin form is then
\begin{eqnarray*}
\delta\varepsilon
& = & \left( \pd{\varepsilon_i}{q^j}
- \half \dT \, \pd{\varepsilon_i}{\dot{q}^j} 
+ \quarter \dT^2 \, \pd{\varepsilon_i}{\ddot{q}^j} \right) dq^j \wedge dq^i \\
& & + \, \half \left( \pd{\varepsilon_i}{\dot{q}^j} 
+ \pd{\varepsilon_j}{\dot{q}^i}
- \dT \pd{\varepsilon_i}{\ddot{q}^j}
- \dT \pd{\varepsilon_j}{\ddot{q}^i} \right) d\dot{q}^j \wedge dq^i \\
& & + \, \quarter \left( \pd{\varepsilon_i}{\ddot{q}^j}
- \pd{\varepsilon_j}{\ddot{q}^i} \right) d\ddot{q}^j \wedge dq^i 
- \half \pd{\varepsilon_i}{\ddot{q}^j}
d\dot{q}^j \wedge d\dot{q}^i \, ,
\end{eqnarray*} 
and so the vanishing of these coefficients is a necessary and sufficient condition for the form to be locally variational. If this is the case we can use the procedure from the previous section to construct a suitable Lagrangian; we finish by demonstrating that making a suitable choice when integrating the two-form $Pd\varepsilon$ allows us to find a Lagrangian which is first-order.

We start by using the condition $\delta\varepsilon = 0$ to tell us about the structure of the functions $\varepsilon_i$. From the term in $d\ddot{q}^j \wedge dq^i$ we have
\[
\pd{\varepsilon_i}{\ddot{q}^j} = \pd{\varepsilon_j}{\ddot{q}^i} \, ;
\]
and then from the term in $d\dot{q}^j \wedge dq^i$ we have just a single term involving $q_{(3)}^k$,
\[
q_{(3)}^k \pdb{\varepsilon_i}{\ddot{q}^k}{\ddot{q}^j}
\]
so we see that $\varepsilon_i$ is affine in the second derivative coordinates. Put
\[
\varepsilon_i = A_{ij} \ddot{q}^j + B_i
\]
where $A_{ij}$ and $B_i$ depend only on $q^k, \dot{q}^k$ and where, from above, $A_{ij} = A_{ji}$.

We now consider the construction of a Lagrangian. First, we see that
\begin{eqnarray*}
Pd\varepsilon & = & \half \pd{\varepsilon_i}{\dot{q}^j} dq^j \wedge dq^i
+ \pd{\varepsilon_i}{\ddot{q}^j} d\dot{q}^j \wedge dq^i \\
& = & \half \left( \ddot{q}^k \pd{A_{ik}}{\dot{q}^j} + \pd{B_i}{\dot{q}^j} \right) dq^j \wedge dq^i
+ A_{ij} d\dot{q}^j \wedge dq^i \, ,
\end{eqnarray*}
and we know that this second-order 2-form must be $d$-closed. Put $Pd\varepsilon = d\kappa$, where in general
\[
\kappa = f_i dq^i + g_i d\dot{q}^i + h_i \ddot{q}^i
\]
will be a second-order 1-form, so that
\begin{eqnarray*}
d\kappa & = & \pd{f_i}{q^j} dq^j \wedge dq^i
+ \pd{f_i}{\dot{q}^j} d\dot{q}^j \wedge dq^i
+ \pd{f_i}{\ddot{q}^j} d\ddot{q}^j \wedge dq^i \\
& & + \pd{g_i}{q^j} dq^j \wedge d\dot{q}^i
+ \pd{g_i}{\dot{q}^j} d\dot{q}^j \wedge d\dot{q}^i
+ \pd{g_i}{\ddot{q}^j} d\ddot{q}^j \wedge d\dot{q}^i \\
& & + \pd{h_i}{q^j} dq^j \wedge d\ddot{q}^i
+ \pd{h_i}{\dot{q}^j} d\dot{q}^j \wedge d\ddot{q}^i
+ \pd{h_i}{\ddot{q}^j} d\ddot{q}^j \wedge d\ddot{q}^i \, .
\end{eqnarray*}
Of course we are free to add to $\kappa$ any exact 1-form $d\mu$, where $\mu$ is a function. But from $d\kappa = Pd\varepsilon$ we have
\[
\pd{h_i}{\ddot{q}^j} - \pd{h_j}{\ddot{q}^i} = 0 \, ,
\]
so we see that there is a function $h$ such that
\[
\quad h_i = \pd{h}{\ddot{q}^i} \, .
\]
We may now consider $\tilde{\kappa} = \kappa - dh$, and again $Pd\varepsilon = d\tilde{\kappa}$. We have
\begin{eqnarray*}
\tilde{\kappa} & = & \left( f_i dq^i + \pd{g}{\dot{q}^i} d\dot{q}^i + \pd{h}{\ddot{q}^i} d\ddot{q}^i \right)
- \left( \pd{h}{q^i} dq^i + \pd{h}{\dot{q}^i} d\dot{q}^i + \pd{h}{\ddot{q}^i} d\ddot{q}^i \right) \\
& = & \left( f_i - \pd{h}{q^i} \right) dq^i
+ \left( \pd{g}{\dot{q}^i} - \pd{h}{\dot{q}^i} \right) d\dot{q}^i \\
& = & \tilde{f}_i dq^i + \tilde{g}_i d\dot{q}^i 
\end{eqnarray*}
so that
\begin{eqnarray*}
d\tilde{\kappa} & = & \pd{\tilde{f}_i}{q^j} dq^j \wedge dq^i
+ \pd{\tilde{f}_i}{\dot{q}^j} d\dot{q}^j \wedge dq^i
+ \pd{\tilde{f}_i}{\ddot{q}^j} d\ddot{q}^j \wedge dq^i \\
& & + \pd{\tilde{g}_i}{q^j} dq^j \wedge d\dot{q}^i
+ \pd{\tilde{g}_i}{\dot{q}^j} d\dot{q}^j \wedge d\dot{q}^i
+ \pd{\tilde{g}_i}{\ddot{q}^j} d\ddot{q}^j \wedge d\dot{q}^i \, .
\end{eqnarray*}
Using the structure of $Pd\varepsilon$ again, we see that $\tilde{f}_i$ and $\tilde{g}_i$ are first-order functions. Furthermore, from
\[
\pd{\tilde{g}_i}{\dot{q}^j} - \pd{\tilde{g}_j}{\dot{q}^i} = 0
\]
we see that there is a function $\tilde{g}$ such that
\[
\tilde{g}_i = \pd{\tilde{g}}{\dot{q}^i} \, ,
\]
and we may now consider $\hat{\kappa} = \tilde{\kappa} - d\tilde{g}$ with, once again, $Pd\varepsilon = d\hat{\kappa}$. We now have
\begin{eqnarray*}
\hat{\kappa} & = & \left( \tilde{f}_i dq^i + \pd{\tilde{g}}{\dot{q}_i} d\dot{q}^i \right)
- \left( \pd{\tilde{g}}{q^i} dq^i + \pd{\tilde{g}}{\dot{q}^i} d\dot{q}^i \right) \\
& = & \left( \tilde{f}_i - \pd{\tilde{g}}{q^i} \right) dq^i \\
& = & \hat{f}_i dq^i
\end{eqnarray*}
so that
\[
d\hat{\kappa} = \pd{\hat{f}_i}{q^j} dq^j \wedge dq^i + \pd{\hat{f}_i}{\dot{q}^j} d\dot{q}^j \wedge dq^i \, .
\]
Thus
\[
\pd{\hat{f}_i}{\dot{q}^j} = A_{ij}
\]
and as $A_{ij} = A_{ji}$ it follows that there is a first-order function $\hat{f}$ such that
\[
\hat{f}_i = \pd{\hat{f}}{\dot{q}^i}
\]
and therefore that
\[
A_{ij} = \pdb{\hat{f}}{\dot{q}^i}{\dot{q}^j} \, .
\]

We now have a suitable 1-form $\hat{\kappa}$, and we know from the theory that $\varepsilon - \dT \hat{\kappa}$ is $d$-closed, so we can find a function $L$ such that $dL = \varepsilon - \dT \hat{\kappa}$. Now
\begin{eqnarray*}
\dT \hat{\kappa} & = & \dT \left( \pd{\hat{f}}{\dot{q}^i} \right) dq^i + \pd{\hat{f}}{\dot{q}^i} d\dot{q}^i \\
& = & \left( \dot{q}^j \pdb{\hat{f}}{q^j}{\dot{q}^i} + \ddot{q}^j \pdb{\hat{f}}{\dot{q}^j}{\dot{q}^i} \right)
dq^i + \pd{\hat{f}}{\dot{q}^i} d\dot{q}^i
\end{eqnarray*}
so that
\begin{eqnarray*}
\varepsilon - \dT \hat{\kappa} & = & \left( \varepsilon_i - \dot{q}^j \pdb{\hat{f}}{q^j}{\dot{q}^i} - \ddot{q}^j \pdb{\hat{f}}{\dot{q}^j}{\dot{q}^i} \right) dq^i - \pd{\hat{f}}{\dot{q}^i} d\dot{q}^i \\
& = & \left( A_{ij} \ddot{q}^j + B_i - \dot{q}^j \pdb{\hat{f}}{q^j}{\dot{q}^i}
- \ddot{q}^j A_{ji} \right) dq^i - \pd{\hat{f}}{\dot{q}^i} d\dot{q}^i \\
& = & \left( B_i - \dot{q}^j \pdb{\hat{f}}{q^j}{\dot{q}^i} \right) dq^i
- \pd{\hat{f}}{\dot{q}^i} d\dot{q}^i \\
& = & \pd{L}{q^i} dq^i + \pd{L}{\dot{q}^i} d\dot{q}^i \, .
\end{eqnarray*}
Thus we see that $L$ is a first-order function, and that
\[
\hat{f} = F - L
\]
for some function $F$ depending only on the coordinates $q^i$. It follows that
\[
\pdb{\hat{f}}{q^j}{\dot{q}^i} = - \pdb{L}{q^j}{\dot{q}^i} \, ,
\]
and so finally we obtain
\begin{eqnarray*}
A_{ij} & = & - \pdb{L}{\dot{q}^i}{\dot{q}^j} \\
B_i & = & \pd{L}{q^i} - \dot{q}^j \pdb{L}{q^j}{\dot{q}^i} \, .
\end{eqnarray*}


\begin{thebibliography}{00}

\bibitem{AndDuch}
I.\ Anderson, T.\ Duchamp:\ On the existence of global variational
principles \textit{Amer.\ J.\ Math.\/} \textbf{102} (1980) 781--867

\bibitem{CS1}
M.\ Crampin, D.J.\ Saunders:\ The Hilbert-Carath\'{e}odory form for
parametric multiple integral problems in the calculus of variations
\textit{Acta Appl.\ Math.\/} \textbf{76(1)} (2003) 37--55

\bibitem{CS2}
M.\ Crampin, D.J.\ Saunders:\ The Hilbert-Carath\'{e}odory and
Poincar\'{e}-Cartan forms for higher-order multiple-integral variational
problems \textit{Houston J.\ Math.\/} \textbf{30(3)} (2004) 657--689

\bibitem{CS3}
M.\ Crampin, D.J.\ Saunders:\ On null Lagrangians \textit{Diff.\ Geom.\ Appl.\/}
\textbf{22(2)} (2005) 131--146

\bibitem{CSC}
M.\ Crampin, W.\ Sarlet, F.\ Cantrijn:\ Higher-order differential equations and higher-order Lagrangian mechanics \textit{Math.\ Proc.\ Camb.\ Phil.\ Soc.\/}
\textbf{99} (2005) 565--587

\bibitem{KrupSeq}
D.\ Krupka: Variational sequences on finite order jet spaces
\textit{In: Differential Geometry and its Applications} (Singapore: World Scientific) 1990

\bibitem{Sau}
D.J.\ Saunders:\ Homogeneous variational complexes and bicomplexes \textit{Preprint} (2005)
ArXiV: math.DG/0512383

\bibitem{Tak}
F.\ Takens:\ A global version of the inverse problem of the calculus of
variations \textit{J.\ Diff.\ Geom.\/} \textbf{14} (1979) 543--562

\bibitem{Tul}
W.M.\ Tulczyjew:\ The Euler-Lagrange Resolution \textit{In: Lecture Notes in Math.
836, Differential Geometric Methods in Mathematical Physics} (Berlin, Springer) 1979

\bibitem{Vin1}
A.M.\ Vinogradov:\ On the algebro-geometric foundations of Lagrangian field
theory \textit{Soviet Math.\ Dokl.\/} \textbf{18} (1977) 1200--1204

\bibitem{Vin2}
A.M.\ Vinogradov:\ A spectral sequence associated with a non-linear
differential equation, and algebro-geometric foundations of Lagrangian field
theory with constraints \textit{Soviet Math.\ Dokl.\/} \textbf{19} (1978) 144--148

\bibitem{Vit}
R.\ Vitolo:\ Finite order Lagrangian bicomplexes \textit{Math.\ Proc.\ Camb.\
Phil.\ Soc.\/} \textbf{125} (1999) 321--333


\end{thebibliography}
\end{document}